\documentclass[a4paper,11pt]{article}
\RequirePackage{amsthm}
\RequirePackage{amssymb}
\usepackage[all,ps]{xy}
\usepackage{latexsym}
\usepackage{amscd}
\newtheorem{thm}{Theorem}[section]
\newtheorem{prop}[thm]{Proposition}
\newtheorem{lem}[thm]{Lemma}

\newtheorem{rmk}[thm]{Remark}

\makeatletter
 \def\ps@notecras{%
    \def\@oddfoot{\small\@date\hfil\slshape\@pauthor\hfil\upshape\thepage}
    \def\@oddfoot{\basdepage}%
   \def\@oddhead{\tetedepage}
  \let\@mkboth\@gobbletwo
 \let\sectionmark\@gobble
\let\subsectionmark\@gobble}
\makeatother
\title{\bf Rational curves on general projective hypersurfaces}
\author{Gianluca Pacienza}
\date{}
\begin{document}
\maketitle
%
%
\begin{quote}
{\small{\bf Abstract.} {In this article, we study the geometry of
$k$-dimensional subvarieties with geometric genus zero of a 
general projective hypersurface $X_d\subset \mathbf P^n$ of degree 
$d=2n-2-k$, where $k$ is an integer such that $1\leq k\leq n-5$. 
As a corollary of our main result
we obtain that the only rational curves lying on the
general hypersuface $X_{2n-3}\subset \mathbf P^n$, for $n\geq 6,$ 
are the lines.}}
\end{quote}

%
%
\section {Introduction}
%

It was shown by H. Clemens [C] that the general (in the countable Zariski
topology) hypersurface of degree $d$ in $\mathbf P^n$ does not contain any
rational curve, if $d$ is sufficiently large.
More precisely, he proved the following:
\bigskip
\\
{\bf Theorem (Clemens).}
 {\it There is no rational curve on the general hypersurface 
 $X_d\subset \mathbf P^n$ of degree $d\geq 2n-1,\ n\geq 3$.}
\bigskip

L. Ein studied more generally (see [E1] and [E2]) 
the geometric genus of subvarieties
contained in complete intersections $X_{(d_1,\dots ,d_r)}\subset M$, where $M$ 
is an arbitrary smooth and projective ambient variety. He proved 
that if $d_1+\dots +d_r\geq 2n-r-k+1$, any $k$-dimensional subvariety
$Y\subset X$ has desingularization with effective canonical bundle.
However, in the case of a
hypersurface $X_d\subset \mathbf P^n$, he obtained
the lower bound $d\geq 2n-k$ on the degree of $X$, which was the same
as Clemens' for $k=1$, and was not optimal. 
Indeed, it was classically known that the lines lying on 
the general hypersurface $X_d\subset \mathbf P^n$ of degree
$d=2n-2-k$ cover a $k$-dimensional subvariety, which
then has geometric genus zero.
Thus nothing was known yet
about the canonical bundle of subvarieties of dimension $k$ on the 
general $X_{2n-1-k}\subset \mathbf P^n$. Voisin ([V2],
[V3]) showed then that it was possible to sharpen Ein's bound 
by one, as conjectured by Clemens himself, by proving:
\bigskip
\\
{\bf Theorem (Voisin).}
 {\it Let $X_d\subset \mathbf P^n$ be a general hypersurface of degree 
 $d\geq 2n-1-k$, where $k$ is an integer such that $1\leq k\leq n-3$.
 Then any $k$-dimensional subvariety $Y$ of $X$ has desingularization
 $\tilde Y $ with effective canonical bundle.}
\bigskip

For $k=1$ we immediately obtain that the general 
$X_{2n-2}\subset \mathbf P^n,\ n\geq 4$, contains no rational curves.
Taking $k=2$ we get another very interesting corollary whose
analogue in the case $n=4$ would solve
Clemens' conjecture on the finiteness of rational curves of fixed 
degree on the general quintic threefold in ${\mathbf P}^4$:
\bigskip
\\ 
{\bf Corollary (Voisin).}
{\it For each integer $\delta\geq 1$, the general hypersurface  
 $X_{2n-3}\subset \mathbf P^n,\ n\geq 5$, contains at most a finite
number of rational curves of degree $\delta$.}
\bigskip

The goal of our work is to investigate, for the general 
$X_{2n-2-k}\subset \mathbf P^n$, \\ $1\leq k\leq n-4$, the geometry of its
$k$-dimensional subvarieties having geometric genus zero.
Since the locus covered by the lines of $X$ is the only known example of
such a subvariety, it seems natural to start with the following: 
\bigskip
\\
{\bf Question.} 
{\it Is the variety covered by the lines the only subvariety 
 of dimension $k$ with geometric genus
 zero on the general hypersurface 
 $X_{2n-2-k}\subset \mathbf P^n,\ 1\leq k\leq n-4$?} 
\bigskip

Remark that the numerical hypothesis $1\leq k\leq n-4$ implies the positivity 
of the canonical bundle of $X_{2n-2-k}\subset \mathbf P^n$, and
gives meaning to the question in contrast to the case of the
Calabi-Yau hypersurface $X_{n+1}\subset \mathbf P^n$.
The main result of this paper gives actually a positive answer to the previous
question for $1\leq k\leq n-5$:
\bigskip
\\
{\bf Theorem.}
 {\it Let $X_d\subset \mathbf P^n$ be a general hypersurface of degree 
 $d=2n-2-k$, where $k$ is an integer such that $1\leq k\leq n-5$.
 Then any subvariety $Y\subset X$ of dimension $k$, whose
 desingularization $\tilde Y$ has $h^0 (\tilde Y,K_{\tilde Y})=0$,
 is a component of the ($k$-dimensional) 
 subvariety covered by the lines lying on $X$.}
\bigskip

Taking $k=1$  we get, for $n\geq 6$,  
a corollary on rational curves on the 
general projective hypersurface
of degree $d=2n-3$. Voisin's corollary already implies that, 
for each fixed integer $\delta \geq 1$,
there are only a finite number of such curves of degree $\delta$.  
Here we prove that there are only lines - whose
number is easily computed as the top Chern class of a
certain vector bundle on the Grassmannian of lines in $\mathbf P^n$:    
\bigskip
\\ 
{\bf Corollary.}
{\it There is no rational curve of degree
$\delta\geq 2$ on the general hypersurface 
$X_{2n-3}\subset \mathbf P^n,\ n\geq 6$.}
\bigskip

Throughout this paper we work on the field of complex number $\mathbb C$.
\bigskip
\\
{\small {\it Acknowledgements.}\ This paper, which is a part of my 
Ph.D. thesis, has greatly benefited from the generous help 
and guidance of my advisor, Prof. Claire Voisin. 
I would also like to thank Prof. A. F. Lopez for his 
encouragement and many useful discussions.
I would like to sincerely thank the referee. 
He or she did a remarkable amount of
work, simplifying some of my proofs, checking any detail,
and suggesting how to put in greater evidence the key 
underlying ideas driving the result. Finally, with great
generosity, he or she shared with us some new ideas of his/hers 
to sharpen the results presented here.}

%
\section {Preliminaries and sketch of the proof}
%

To motivate our approach  and render our proof
more transparent, we will briefly sketch 
the key ideas contained in [E1], [V2] and [V3].
We start with some 
\bigskip
\\
{\bf Notation.}

\noindent $S^d:=H^0 (\mathbf P^n, {\cal  O}_{\mathbf P^n}(d))$;
\\ 
$S^d_x:=H^0 (\mathbf P^n, {\cal I}_x\otimes{\cal  O}_{\mathbf P^n}(d))$; 
\\ 
$N:=h^0(\mathbf P^n, {\cal  O}_{\mathbf P^n}(d))= dim\ S^d$;
\\ 
${\cal X}\subset {\mathbf P^n}\times S^d$ 
will denote the universal hypersurface of degree $d$;
\\
$X_F\subset {\mathbf P^n}$ the fiber of the family ${\cal X}$
over $F\in S^d$, i.e. the hypersurface 
defined by $F$. 
\smallskip

Let $U\rightarrow S^d$ be an \'etale map and
${\cal Y}\subset {\cal X}_U$ a universal, reduced and 
irreducible subscheme of relative dimension $k$
(in the following, by abuse of notation, we will often omit
the \'etale base change).
We may obviously assume $\cal Y$ invariant under some 
lift of the natural action of $GL(n+1)$ on  ${\mathbf P^n}\times S^d$:
$
g(x,F)= (g(x), (g^{-1})^*F)),\ g\in GL(n+1).
$
Let $\tilde {\cal Y}\rightarrow {\cal Y}$ be a desingularization and 
$\tilde {\cal Y}{\buildrel j\over\longrightarrow}{\cal X}_U$ 
the natural induced map.
Let $\pi:{\cal X}\rightarrow {\mathbf P^n}$ be the projection on the
first component and  $T^{vert}_{\cal X}$ (resp. $T^{vert}_{\cal Y}$)
the vertical part of $T_{\cal X}$ (resp. $T_{\cal Y}$) 
w.r.t. $\pi$, i.e. $T^{vert}_{\cal X}$ (resp. $T^{vert}_{\cal Y}$) is the sheaf
defined by
$$
 0\rightarrow T^{vert}_{\cal X}\rightarrow T{\cal X}{\buildrel {\pi_*}\over
 \longrightarrow} T{\mathbf P^n}\rightarrow 0
$$
$$
 ({\it {resp.}}\ \ \   0\rightarrow T^{vert}_{\cal Y}\rightarrow 
 T{\cal Y}\ {\buildrel {\pi_*}\over \longrightarrow} T{\mathbf P^n}).
$$
The hypothesis of $GL(n+1)$-invariance for $\cal Y$ 
has two easy
but very important consequences that  
will be  frequently used in the rest of the paper:
\begin{lem}\label{GL}
Let ${\cal Y}$ be a subvariety of ${\cal X}\times {\mathbf P^n}$
of relative dimension $k$ and invariant under the action of $GL(n+1)$.
Then:

(i) 
$
 codim_{\ T^{vert}_{{\cal X},(y,F)}}\ T^{vert}_{{\cal Y},(y,F)}
 =codim_{\cal X}{\cal Y}=n-k-1;
$
in particular, since we are assuming $1\leq k\leq n-5$, we have that
$$
 codim\ _{\cal X}\ {\cal Y}\ \geq 4.
$$

(ii) 
$
 T^{vert}_{{\cal Y},(y,F)}\supset\ <S^1_y \cdot J_F^{d-1}, F>,
$
where $J_F^{d-1}$ is the Jacobian ideal of $F$.
\end{lem}
\begin{proof}
(i) Use the surjectivity of the map 
$T{\cal Y}\ {\buildrel {\pi_*}\over 
\longrightarrow} T{\mathbf P^n}$.

\noindent
(ii) It follows from the fact that, by $GL(n+1)$-invariance,
$T^{vert}_{{\cal Y},(y,F)}$ contains the vertical part of the
tangent space to the orbit 
of $(y,F)$ under the action of $GL(n+1)$.
\end{proof}

Assume $d\geq 2n-k$ and let
$Y_F\subset X_F$ be a general fiber of the subfamily 
${\cal Y}\subset {\cal X}_U$. 
Then to prove Clemens' result (which corresponds to
the case $k=1$)
we have to show that  
$h^0 (\tilde Y_F, K_{\tilde Y_F})\not =0$, where 
$\tilde Y_F\rightarrow Y_F$ is the desingularization. 

Ein's idea is to produce, by adjunction, a non-zero section in 
$H^0 (\tilde Y_F, K_{{\tilde Y}_F})$ by restricting to $\tilde {\cal Y}$ 
holomorphic forms on ${\cal X}_U$ defined along $X_F$. 
The main technical difficulty consists in controlling the positivity
of the twisted ``vertical'' tangent bundle to the universal 
hypersurfaces. We start then recording, in the first subsection, 
the needed positivity result, 
and an equivalent one for a bundle on the Grassmannian 
of lines in ${\mathbf P}^n$ that will be used
later.

%
\subsection {Positivity results}
%
Let $d$ be a positive integer.
Consider the bundle $M^d_{{\mathbf P}^n}$ defined by the exact sequence
\begin{eqnarray}\label{M_d}
 0\rightarrow M_{{\mathbf P}^n}^d
 \rightarrow S^d \otimes {\cal O}_{\mathbf P^n}
 {\buildrel {ev}\over
 \longrightarrow} {\cal O}_{\mathbf P^n}(d)\rightarrow 0,
\end{eqnarray}
whose fiber at a point $x$ identifies by definition with $S_x^d$. 
From the inclusion ${\cal X}\hookrightarrow {\mathbf P^n}\times S^d$ 
we get the exact sequence
$$
 0\rightarrow T{\cal X}_{|X_F}\rightarrow 
 T{\mathbf P^n}_{|X_F}\oplus (S^d\otimes {\cal O}_{X_F})
 \rightarrow {\cal O}_{X_F}(d)\rightarrow 0,
$$
which combined with (\ref{M_d}) gives us
$$
 0\rightarrow {M_{{\mathbf P}^n}^d}_{|X_F}\rightarrow T{\cal X}_{|X_F}
 \rightarrow T{\mathbf P^n}_{|X_F}\rightarrow 0.
$$
In other words 
${M_{{\mathbf P}^n}^d}_{|X_F}$ identifies to the 
vertical part of $T_{\cal X}\otimes {\cal O}_{X_F}$ with respect to the 
projection to ${\mathbf P^n}$.
 
Let $G:=Grass(1,n)$ be the Grassmannian of lines in ${\mathbf P}^n$,
${\cal O}_G (1)$ the line bundle on $G$ giving its Pl\"ucker polarization,
and ${\cal E}_d$ be the $d^{th}$-symmetric power of the dual of the 
tautological subbundle on $G$. 
Recall that the fibre of ${\cal E}_d$ at a point $[{\ell}]$ 
is, by definition, given by $H^0 ({\ell},{\cal O}_{\ell} (d))$.

Let $M_G^d$ be the vector bundle on $G$ 
defined as the kernel of the evaluation map:
$$
 0\rightarrow M_G^d\rightarrow S^d \otimes {\cal O}_G 
 \rightarrow {\cal E}_d \rightarrow 0.
$$
Notice that the fiber of $M_G^d$ at a point 
$[\ell]$ is equal to $I_{\ell}(d):=H^0 ({\cal I}_{\ell}(d))$. 

Then we have the following 

\begin{prop}\label{pos}\footnote{The quick proof of (i) we reproduce 
here is due to the referee, who pointed out to us the possibility
of using Bott's and Mumford's results. The referee also suggested us
the proof of (ii), which is much simpler than the one we originally proposed.}
 (i) 
 $M_{{\mathbf P}^n}^d \otimes {\cal O}_{{\mathbf P}^n}(1)$ is 
 generated be its global sections;\\
 (ii) $M_G^d\otimes {\cal O}_G (1)$ is 
 generated be its global sections.
\end{prop}
\begin{proof}
(i) We start observing that the sheaf 
$\Omega^s_{{\mathbf P}^n} (s+1)$ is globally
generated. Indeed, by Bott's vanishing theorem [B]
we have
\begin{equation}
 H^i (\Omega^p_{{\mathbf P}^n}(q))=0,\ \forall i>0\ q\geq p+1-i.
\end{equation}  
Therefore, by Mumford's m-regularity theorem ([M1], page 99) the maps
$$
 H^{0}(\Omega _{\mathbf{P}^{n}}^{p}(q)) \otimes
 H^{0}(\mathcal{O}_{\mathbf{P}^{n}}(1)) \rightarrow
 H^{0}(\Omega _{\mathbf{P}^{n}}^{p}(q+1)),
$$
are surjective for $q\geq p+1$. The result now follows immediately
recalling the isomorphism
$$
 M_{{\mathbf P}^n}^1\cong\Omega _{\mathbf{P}^{n}}^{1}(1)
$$  
and the surjection
$$
 S^{d-1}\otimes M_{{\mathbf P}^n}^1\twoheadrightarrow M_{{\mathbf P}^n}^d.
$$

(ii) Again it suffices to check the case $d=1$. Using the irreducibility of
the action of $GL\left( n+1\right)$ on the Grassmannian,
 it suffices to construct a single
non-trivial meromorphic section of $M_{G}^{1}$ with simple pole along the
zero set of a Pl\"{u}cker coordinate. To do this, for all lines $\ell$ not
meeting
$$
X_{1}=X_{2}=0
$$
the Pl\"{u}cker coordinate $p_{12}(\ell) \neq 0$ so there is, by
Cramer's rule, a unique $L=X_{0}+aX_{1}+bX_{2}\in M_{G}^{1}$ containing $\ell$.

\end{proof}

%
\subsection {Proof of Ein's and Voisin's results}
%
%
%
%
%

Following Ein [E1], [E2], one can use the
positivity result (\ref{pos}), (i), to produce holomorphic forms on the 
(vertical) tangent space to the family ${\cal X}_U$. 
Then, by pulling back them to
${\tilde Y_F}$ and using the adjunction formula, it will be possible to
provide a non zero section
of $H^0 (\tilde Y_F, K_{{\tilde Y}_F})$.
To make this more precise, first
recall the following elementary facts:
\smallskip

(i) ${\Omega ^{N+k}_{{\tilde{\cal Y}}}}_{|{\tilde Y_F}}
\cong K_{\tilde Y_F}$;
\smallskip

(ii) $(\wedge^{n-1-k} T{{\cal X}_U}_{|{X_F}})\otimes K_{X_F}
\cong {\Omega ^{N+k}_{{\cal X}_U}}_{|{X_F}}$.
\smallskip

\noindent Therefore, from the natural morphism 
$\Omega ^1_{{{\cal X}_U}}\rightarrow \Omega ^1_{{\tilde{\cal Y}}}$, we
get a map
\begin{equation}\label{restriction}
 (\wedge^{n-1-k} T{{\cal X}_U}_{|{X_F}})\otimes K_{X_F}
 \cong {\Omega ^{N+k}_{{\cal X}_U}}_{|{X_F}}
 \rightarrow {\Omega ^{N+k}_{{\tilde{\cal Y}}}}_{|{\tilde Y_F}}
 \cong K_{\tilde Y_F}.
\end{equation}
Since $K_{X_F}={\cal O}_{X_F}(d-n-1)={\cal O}_{X_F}((n-k-1)+(d-2n+k))$ 
and
$$
\wedge^{n-1-k} T{{\cal X}_U}_{|{X_F}}(n-k-1)=
\wedge^{n-1-k} (T{{\cal X}_U}_{|{X_F}}(1)),
$$
we have
\begin{equation}\label{wedge}
 (\wedge^{n-1-k} T{{\cal X}_U}_{|{X_F}})\otimes K_{X_F}=
 \wedge^{n-1-k} (T{{\cal X}_U}_{|{X_F}}(1))\otimes 
 {\cal O}_{X_F}(d-2n+k).
\end{equation}

Now, since we are supposing $d\geq 2n-k$,
Proposition (\ref{pos}) (i) implies that the vertical part of
$$
 \wedge^{n-1-k} (T{{\cal X}_U}_{|{X_F}}(1))\otimes {\cal O}_{X_F}(d-2n+k)
 \cong {\Omega ^{N+k}_{{\cal X}_U}}_{|{X_F}},
$$
namely, the subsheaf
$$
 \wedge^{n-1-k} ({M_{{\mathbf P}^n}^d}_{|{X_F}})\otimes K_{X_F}
 =\wedge^{n-1-k} ({M_{{\mathbf P}^n}^d}_{|{X_F}}(1))
 \otimes {\cal O}_{X_F}(d-2n+k),
$$
is globally generated. 
Composing the inclusion 
$$
 \wedge^{n-1-k} ({M_{{\mathbf P}^n}^d}_{|{X_F}})\otimes K_{X_F}
 \hookrightarrow
 \wedge^{n-1-k} (T{{\cal X}_U}_{|{X_F}})\otimes K_{X_F}
$$ 
with the restriction map defined in (\ref{restriction}), we have a 
natural morphism
\begin{equation}\label{res}
 \wedge^{n-1-k} ({M_{{\mathbf P}^n}^d}_{|{X_F}})\otimes K_{X_F}
 \to K_{\tilde Y_F}.
\end{equation}

Ein's result is then given by the following

\begin{lem}
 Let $F$ be a general polynomial
 of degree $2n-k$. The map
 $$
   H^0 (\wedge^{n-1-k} ({M_{{\mathbf P}^n}^d}_{|{X_F}})\otimes K_{X_F})
   \to H^0 ({K_{\tilde Y_F}}),
 $$
 induced in cohomology by (\ref{res}), is non zero. 
\end{lem}
\begin{proof}
 By Lemma (\ref{GL}), (i), we have 
 $$
  codim_{\ T^{vert}_{{\cal X},(y,F)}}\ T^{vert}_{{\cal Y},(y,F)}
  =codim_{\cal X}{\cal Y}.
 $$
 Let $(y,F)$ be a smooth point of ${\cal Y}$. 
 Since the bundle
 $\wedge^{n-1-k} ({M_{{\mathbf P}^n}^d}_{|{X_F}})\otimes K_{X_F}$
 is generated by its global sections, there exists a section
 $$
  s\in 
  H^0 (\wedge^{n-1-k} ({M_{{\mathbf P}^n}^d}_{|{X_F}})\otimes K_{X_F})
 $$ 
 such that
 \begin{equation}
  <s(y), T^{vert}_{{\cal Y},(y,F)}> \not= 0.
 \end{equation}
 Since $j:\tilde {\cal Y}\rightarrow {\cal X}_U$ is
 generically an immersion, we obtain from the above a non zero
 element in  $H^0 ({K_{\tilde Y_F}})$ coming from 
 $H^0 (\wedge^{n-1-k} ({M_{{\mathbf P}^n}^d}_{|{X_F}})\otimes K_{X_F})$.
\end{proof}

(For other proofs of Clemens' theorem see, of course, [C] 
and also [CLR]).

%
%
%
%

In order to try and improve the bound on the degree by one, we
observe that, if $d=2n-1-k$, then 
$
K_{X_F}={\cal O}_{X_F}(n-2-k),
$
so we have, as in (\ref{restriction}), a map
\begin{equation}\label{shift}
 \wedge^{n-1-k} T{{\cal X}_U}_{|{X_F}}(n-k-2)
 \cong {\Omega ^{N+k}_{{\cal X}_U}}_{|{X_F}}\rightarrow
 {\Omega ^{N+k}_{{\tilde{\cal Y}}}}_{|{\tilde Y_F}}
 \cong K_{\tilde Y_F}.
\end{equation}
As we saw in Lemma 2.1, (i), 
by the hypothesis of $GL(n+1)$-invariance on ${\cal Y}$, 
the relevant part of the tangent space to look at is the vertical one,
hence we focus our attention on the map
\begin{equation}\label{vertshift}
\wedge^{n-1-k}{M_{{\mathbf P}^n}^d}_{|{X_F}}(n-k-2)
 \rightarrow
 {\Omega ^{N+k}_{{\tilde{\cal Y}}}}_{|{\tilde Y_F}}
 \cong K_{\tilde Y_F}.
\end{equation}

Now, because of the shift between the exterior power and the 
degree of the canonical bundle we are tensoring with on the lefthand
side of (\ref{vertshift}), the global generation of the sheaf 
$\wedge^{n-1-k}{M_{{\mathbf P}^n}^d}_{|{X_F}}(n-k-2)$ 
will not follow from 
the global generation of ${M_{{\mathbf P}^n}^d}(1)$. Voisin's idea is
then to study the positivity of
$H^0 (\wedge^{2}{M_{{\mathbf P}^n}^d} (1))$, to produce
holomorphic forms on the (vertical) tangent space to the
universal hypersurface,
and use the commutative diagram below to produce 
sections in $H^0 (K_{\tilde Y_F})$:
\begin{eqnarray}\label{composite}
 &\ \ \ \ \ \ \ \ \ \ \ \ \ \ \ \ \ \ \ 
 H^0 (\wedge^{n-1-k}{M_{{\mathbf P}^n}^d} _{|{X_F}}(n-2-k))&
 \longrightarrow   \ \ \ \ \ \ \ 
 H^0 (K_{\tilde Y_F})\nonumber\\
 &\uparrow&\ \ \nearrow\nonumber\\
 &H^0 (\wedge^{n-3-k}{M_{{\mathbf P}^n}^d} _{|{X_F}}(n-3-k))\otimes 
 H^0 (\wedge^{2}{M_{{\mathbf P}^n}^d} _{|{X_F}}(1))&
\end{eqnarray}
(the vertical map in (\ref{composite})
is simply obtained by wedging the sections of the sheaves 
$\wedge^{n-3-k} {M_{{\mathbf P}^n}^d}_{|{X_F}}(n-3-k)$ 
and $\wedge^{2} {M_{{\mathbf P}^n}^d}_{|{X_F}}(1)$).
Unfortunately the following fact holds:
\smallskip
\\
{\bf Fact (Amerik-Voisin).} 
{\it $\wedge^{2} {M_{{\mathbf P}^n}^d}(1)$ 
 is not generated by its global sections.}
\smallskip

Indeed, in [V3] the following counterexample to the global generation 
of $\wedge^{2} {M_{{\mathbf P}^n}^d}(1)$ is given. Consider the subvariety
\begin{eqnarray}\label{S_F}
 \Delta_{d,F} :=\lbrace x\in X_F : there\ exists\ a\ line\ {\ell}
 \ s.t.\ {\ell}\cap X_F= d\cdot x\rbrace.
\end{eqnarray}
An elementary dimension count shows that, for generic $F$,
$$
 dim\ \Delta_{d,F}=2n-2-(d-1)=2n-2-(2n-1-k-1)=k
$$ 
(these subvarieties
are generically empty for $d\geq 2n-1$, which is the reason they don't come
into play in Clemens' and Ein's case).  Let  ${\Delta}_d$ be the family
of the $\Delta_{d,F}$'s, let 
$\tilde {\Delta}_d\rightarrow {\Delta}_d$ be a desingularization,
and $j:\tilde {\Delta}_d\rightarrow {\cal X}$ the natural morphism. 
Notice that $\Delta_{d,F}$ parametrizes $0$-cycles of $X_F$ which are all
rationally equivalent 
since, by definition, 
$d\cdot x\equiv H^{n-1}.X_F$, $\forall x\in\Delta_{d,F}$,
where $H$ is the hyperplane divisor in $\mathbf P^n$.
Thus, the variational (and higher dimensional) version of Mumford's
fundamental result on $0$-cycles on surfaces applies (see [M2], and [V1]
for the variational generalization in dimension $2$), 
so we have 
$$
 j^*s=0\ \ in\ \ H^0 ({\Omega ^{N+k}_{{\tilde{\cal Y}}}}_{|{\tilde Y_F}}),\ 
 \forall s\in H^0 ({\Omega ^{N+k}_{\cal X}}_{|{X_F}}),
$$
i.e. the map
$$
 H^0 (\wedge^{n-1-k} T{{\cal X}}_{|{X_F}}(n-2-k))\cong
 H^0 ({\Omega ^{N+k}_{\cal X}}_{|{X_F}})\rightarrow
 H^0 ({\Omega ^{N+k}_{{\tilde{\cal Y}}}}_{|{\tilde Y_F}}) 
 \cong H^0 (K_{\tilde Y_F})
$$
is identically zero and then so is
$$
 H^0 (\wedge^{n-1-k} {M_{{\mathbf P}^n}^d}_{|{X_F}}(n-2-k))
 \to H^0 (K_{\tilde Y_F}).
$$
In particular, by (\ref{composite}) and Proposition \ref{pos}, (i),
we have that, at a smooth point
$(y,F)\in \cal Y$, all the global sections of the bundle 
$\wedge^2 {M_{{\mathbf P}^n}^d}(1)_{|_{X_F}}$, seen as
a line bundle on the Grassmannian of
codimension two subspaces of $T{\cal X}^{vert}_{|{X_F}}$,
vanish on the codimension two subspaces of $T^{vert}_{{\cal X},(y,F)}$ 
containing $T^{vert}_{{\cal Y},(y,F)}$.

Voisin's alternative approach to the problem, as developped in [V3],
consists then in studying the base locus 
of $H^0 (\wedge^2 {M_{{\mathbf P}^n}^d}(1)_{|_{X_F}})$, to investigate
the geometry of 
the subvarieties for which the composite map
in (\ref{composite}) fails to provide non-zero sections of their 
canonical bundle. She shows in [V3] that, 
in the case $d=2n-1-k$, the subvariety $\Delta_{d,F}$
defined in (\ref{S_F}) is the only one for which this phenomenon occurs.
Then, she completes her proof by verifying 
that each component of $\Delta_{d,F}$ has positive geometric genus.

%
\subsection {The strategy of our proof}
%

Our purpose is to study, for $d=2n-2-k$, the geometry of 
$k$-dimensional subvarieties of $X_F\subset \mathbf P^n$,
having geometric genus equal to zero. 
Recall that, since $d=2n-2-k$, we have
$K_{X_F}={\cal O}_{X_F}(n-3-k)$, and note that the composite map
\begin{eqnarray}\label{compositebis}
 &\ \ \ \ \ \ \ \ \ \ \ \ \ \ \ \ \ \ \ 
 H^0 (\wedge^{n-1-k}{M_{{\mathbf P}^n}^d} _{|{X_F}}(n-3-k))&
 \longrightarrow   \ \ \ \ \ \ \ 
 H^0 (K_{\tilde Y_F})\nonumber\\
 &\uparrow&\ \ \nearrow\nonumber\\
 &H^0 (\wedge^{n-5-k}{M_{{\mathbf P}^n}^d} _{|{X_F}}(n-5-k))\otimes 
 H^0 (\wedge^{4}{M_{{\mathbf P}^n}^d} _{|{X_F}}(2))&
\end{eqnarray}

is obviously zero, since we are supposing $h^0 (K_{\tilde Y_F})=0$.
Then the proof of our theorem will naturally be divided into two steps.
In section $3$, analysing the base locus of 
$H^0 (\wedge^{4}{M_{{\mathbf P}^n}^d} _{|{X_F}}(2))$, considered as the
space of sections of a line bundle on the Grassmannian of
codimension four subspaces of $T{\cal X}^{vert}_{|{X_F}}$, we will prove
\bigskip
\\
{\bf Proposition A.} 
{\it Let $X_F\subset \mathbf P^n$ be a general hypersurface of 
 degree $d=2n-2-k$, and $Y_F\subset X_F$ a subvariety
 of dimension $k$ such that $H^0 (\tilde {Y}_F,K_{{\tilde Y}_F})=0$, 
 where $\tilde {Y}_F$ is a 
 desingularization of $Y_F$. Then $Y_F$ has to be contained in
 $$  
  {\Delta}_{d,F}=\lbrace x\in X_F : there\ exists\ a\ line\ \ell\ 
  s.t.\ \ell\cap X_F= d{\cdot} x\rbrace,
 $$
 a subvariety of $X_F$ of dimension $2n-2-(d-1)=k+1$.}
\bigskip

In $\S\ 4$, we will study an
explicit desingularization $\tilde {\Delta}_{d,F}$ 
of ${\Delta}_{d,F}$, given by the zeroes
of a section of a bundle on the incidence variety in 
$\mathbf P^n\times Grass(1,n)$. Denote by $\tilde {\Delta}_{d}$ the
family $({\Delta}_{d,F})_{F\in S^d}$, and recall 
the isomorphism 
$$
 T{\tilde {\Delta}}_{d|\tilde {\Delta}_{d,F}}
 \otimes K_{\tilde {\Delta}_{d,F}}\cong
 {{\Omega}^{N+k}_{\tilde {\Delta}_d}}_{|{{\tilde S}_F}}.
$$
The positivity result (\ref{pos}), (ii), for the bundle  
$M_G^d\otimes {\cal O}_G (1)$ on $Grass(1,n)$, will allow us
to construct a subbundle of 
$
T{\tilde {\Delta}}_{d|\tilde {\Delta}_{d,F}}
\otimes K_{\tilde {\Delta}_{d,F}}
$
generated by its global sections. Using this fact, together with 
the vanishing of the natural restriction map
\begin{equation}\label{S-restriction}
 H^0 (T{\tilde {\Delta}}_{d|\tilde {\Delta}_{d,F}}
 \otimes K_{\tilde {\Delta}_{d,F}})= 
 H^0 ({{\Omega}^{N+k}_{\tilde {\Delta}_d}}_{|{{\tilde S}_F}})
 \rightarrow H^0 (K_{{\tilde Y}_F}),
\end{equation}
we will prove 
\bigskip
\\
{\bf Proposition B.}
{\it Let $F$ be a general polynomial of degree $d=2n-2-k$.
 Let $Y_F\subset {\Delta}_{d,F}$ be a subvariety of codimension $1$ 
 such that $H^0 ({\tilde Y}_F,K_{{\tilde Y}_F})=0$, where 
 ${\tilde Y}_F$ is a desingularization of $Y_F$. 
 Then $Y_F$ has to be a component
 of the $k$-dimensional subvariety of ${\Delta}_{d,F}$ covered by the 
 lines lying on $X_F$.}
\bigskip

These propositions will combine to prove our main theorem.
%
\section 
{Base locus of $\wedge^{4} M_{{\mathbf P}^n}^d (2)$ 
and osculating lines}
%
Let ${\cal X}\subset {\mathbf P^n}\times S^d$ be the universal 
hypersurface of degree $d=2n-2-k$, 
$U\rightarrow S^d$ an \'etale map and ${\cal Y}\subset {\cal X}_U$ 
a universal, reduced and irreducible subscheme of relative dimension
$k$ (to simplify the notation, in what follows we will occasionally omit
the \'etale base change). 
Assume $\cal Y$ invariant under some 
lift of the action of $GL(n+1)$, denote by 
$\tilde {\cal Y}\rightarrow \cal Y$ a desingularization, and suppose that
the fibres of $\tilde {\cal Y}$ verify $h^0 ({\tilde Y_F},K_{\tilde Y_F})=0$.

Consider the bundle $M_{{\mathbf P}^n}^d$ defined by the exact sequence
$$
 0\rightarrow M_{{\mathbf P}^n}^d\rightarrow S^d \otimes {\cal O}_{\mathbf P^n}
 {\buildrel {ev}\over
 \longrightarrow} {\cal O}_{\mathbf P^n}(d)\rightarrow 0,
$$
whose fiber at a point $x$ identifies by definition with $S_x^d$. 
Recall from \S 2.1 that
\begin{equation}\label{vert}
{M_{{\mathbf P}^n}^d}_{|X_F}=T^{vert}_{\cal X}\otimes {\cal O}_{X_F},
\end{equation}
where 
$T^{vert}_{\cal X}$ is the sheaf defined by
$$
 0\rightarrow T^{vert}_{\cal X}\rightarrow T{\cal X}{\buildrel {\pi_*}\over
 \longrightarrow} T{\mathbf P^n}\rightarrow 0.
$$

From the vanishing of the composite map
\begin{eqnarray}\label{compositeter}
 &\ \ \ \ \ \ \ \ \ \ \ \ \ \ \ \ \ \ \ 
 H^0 (\wedge^{n-1-k} T{{\cal X}_U}_{|{X_F}}(n-3-k))&
 \longrightarrow   \ \ \ \ \ \ \ 
 H^0 (K_{\tilde Y_F})\nonumber\\
 &\uparrow&\ \ \nearrow\nonumber\\
 &H^0 (\wedge^{n-5-k} T{{\cal X}_U}_{|{X_F}}(n-5-k))\otimes 
 H^0 (\wedge^{4} T{{\cal X}_U}_{|{X_F}}(2))&
\end{eqnarray}
and (\ref{vert}) we deduce that the composite map
\begin{eqnarray}\label{composite4}
 &\ \ \ \ \ \ \ \ \ \ \ \ \ \ \ \ \ \ \ 
 H^0 (\wedge^{n-1-k}{M_{{\mathbf P}^n}^d} _{|{X_F}}(n-3-k))&
 \longrightarrow   \ \ \ \ \ \ \ 
 H^0 (K_{\tilde Y_F})\nonumber\\
 &\uparrow&\ \ \nearrow\nonumber\\
 &H^0 (\wedge^{n-5-k}{M_{{\mathbf P}^n}^d} _{|{X_F}}(n-5-k))\otimes 
 H^0 (\wedge^{4}{M_{{\mathbf P}^n}^d} _{|{X_F}}(2))&
\end{eqnarray}
is also zero.
Since, by Lemma \ref{pos}, (i),
$\wedge^{n-5-k}{M_{{\mathbf P}^n}^d} _{|{X_F}}(n-5-k)$ is generated
by its global sections, the vanishing of the composite map in 
(\ref{composite4}) and the $GL(n+1)$ invariance of ${\cal Y}$ implies that, 
at a smooth point $(y,F)\in {\cal Y}$, 
any codimension four subspace of $T^{vert}_{{\cal X}_U,(y,F)}=S^d_y$ 
containing $T^{vert}_{{\cal Y},(y,F)}$ is in the base locus
of $H^0 (\wedge^4 {M_{{\mathbf P}^n}^d} _{|{X_F}}(2))$, 
considered as a space of 
sections of a line bundle over the Grassmannian of codimension
four subspaces of $T^{vert}{\cal X}_{|_{X_F}}$. Studying this base locus 
we will see how, at each point $y$ of a subvariety $Y_F\subset X_F$ with
zero geometric genus, the ideal of a line through $y$  
naturally comes into play. More precisely, we will prove
\begin{prop}\label{ideals}
 Let ${\cal Y}\subset {\cal X}_U$ be such that the 
 composite map in (\ref{composite4}) is zero. Then,
 at a smooth point $(y,F)$, the vertical 
 tangent space $T^{vert}_{{\cal Y},(y,F)}$, which is a subspace of
 $T^{vert}_{{\cal X}_U,(y,F)}= S^d_y$, has to contain (at least)
 a hyperplane $H_{\ell_{(y,F)}}\subset I_{\ell _{(y,F)}}(d)$,
 where ${\ell _{(y,F)}}$ is a line passing through $y$.
\end{prop}
We will then study the distribution 
${\cal H}\subset T^{vert}_{{\cal Y}}$, pointwise defined by 
$H_{\ell _{(y,F)}}$, and prove its integrability. The description of the
corresponding foliation and the $GL(n+1)$-invariance of ${\cal Y}$
will allow us to conclude that the line ${\ell}_{(y,F)}$ is such that
$
 {\ell}_{(y,F)}\cap X_F= d\cdot y,
$
thus proving Proposition A.

%
\subsection{Proof of Proposition 3.1}
%
We start with the following
\begin{lem}\label{BL}
 Let $T$ be a codimension four subspace of $S^d_x=(T_{{\cal X},(x,F)})^{vert}$ 
 which is in the base locus of 
 $H^0 (\bigwedge^4 {M_{{\mathbf P}^n}^d}(2))$. 
 Then $T$ has to contain (at least) 
 a hyperplane of $I_{\ell} (d)$, where
 $\ell$ is a line passing through $x$.
\end{lem}
\begin{proof}
 Recall that $H^0 (\bigwedge^2 {M_{{\mathbf P}^n}^d}(1))$ can
 be naturally interpreted as the kernel of the Koszul map 
 $\bigwedge^2 S^d\otimes S^1 \rightarrow S^d \otimes S^{d+1}$. Hence
 one easily verifies that 
 $Im \ (H^0 (\bigwedge^2 {M_{{\mathbf P}^n}^d}(1))
 \rightarrow \bigwedge^2 {M_{{\mathbf P}^n}^d}_{,x})$ 
 contains $PA_1\wedge PA_2$, for all $P\in S^{d-1}$ and $A_i\in S^1_x$.
 Then 
 $$
  Im  \ (H^0 (\bigwedge^4 {M_{{\mathbf P}^n}^d}(2)))
  \rightarrow \bigwedge^4 {M_{{\mathbf P}^n}^d}_{,x})
 $$ 
 contains elements of the form 
 $$
  PA_1\wedge PA_2\wedge QB_1\wedge QB_2,
 $$ 
 for all $P,Q\in S^{d-1}$ and $A_i, B_i\in S^1_x$, 
 coming from the wedge product of elements in 
 $Im \ (H^0 (\bigwedge^2 {M_{{\mathbf P}^n}^d}(1))
 \rightarrow \bigwedge^2 {M_{{\mathbf P}^n}^d}_{,x})$. 
 Since we are supposing that $T$ is in the base locus of 
 $H^0 (\bigwedge^4 {M_{{\mathbf P}^n}^d} (2))$, 
 the previous fact implies in particular that the dimension of the subspace 
 $\lbrace A P: A\in S^1_x\rbrace\ modulo\ T$ is at most $3$,
 i.e. the multiplication map
 \begin{eqnarray}
  m_P : &S^1_x& \rightarrow S^d_x/T \nonumber
  \\&A&\mapsto A\cdot P\ mod\ T\nonumber
 \end{eqnarray}
 cannot be surjective, for any $P\in S^{d-1}$. 
 \\Recall that if $V$ and $W$ are vector spaces, and 
 $Z_k:=\lbrace \phi\in Hom(V,W): rank\ \phi\leq k\rbrace$,
 then 
 \begin{equation}\label{tg}
  T_{{Z_k},\phi}=\lbrace \psi\in Hom(V,W): 
  \psi (ker\phi)\subset Im\phi\rbrace.
 \end{equation}
 If, for generic $P$, the map $m_P$ has rank one, from (\ref{tg})
 we obtain that
 $Q\cdot Ker\ {m_P}\ mod\ T\subset Im\ m_P$, for any $Q\in S^{d-1}$, 
 i.e. $I_{\ell_P} (d)\subset P\cdot S^1_x + T$, where
 ${\ell_P}$ is the line determined by $Ker\ {m_P}$. 
 Then $T$ contains a hyperplane of $I_{\ell _P} (d)$ and the Lemma is proved. 
 
 Thus, we can assume that, for generic $P$. the map $m_P$ has rank 
 at least two.
 Let $A_1, A_2\in S^1_x$ such that 
 $T':= <A_1P, A_2P, T>$ is of codimension $2$ in $S^d_x$. 
 For generic $Q\in S^{d-1}_x$, consider the map
 $$
  m_Q : S^1_x \rightarrow S^d_x/T',
 $$
 whose rank is then equal to $0$ or $1$.
 In the former case $T'$ would then contain 
 $S^{d-1}\cdot Ker\ m_Q=S^{d-1}\cdot S^1_x= S^d_x$, which is absurd since 
 $T'$ has codimension $2$. Hence we can suppose $rk\ m_Q=1$.
 Then by [V2], Lemma $2.3$, 
 $T'$ contains the degree d part of the ideal of a line 
 ${\ell_Q}$ passing through $x$, 
 and hence $T$ contains at least a codimension two subspace of 
 $I_{\ell_Q} (d)$. Assume first that the line does not vary with $Q$, and
 denote it by $\ell$.
 If 
 $$
  codim_{I_{\ell} (d)}\ T\cap I_{\ell} (d)=2,
 $$
 then the image $\overline T$ of $T$ in $H^0 ({\cal O}_{\ell}(d))$
 has codimension 2 in $H^0 ({\cal O}_{\ell}(d)(-x))$.
 On the other hand, since $T'$ contains $I_{\ell} (d)$, its
 reduction $\overline {T'}$ modulo $I_{\ell} (d)$ has also
 codimension 2 in $H^0 ({\cal O}_{\ell}(d)(-x))$. 
 Hence $\overline T=\overline {T'}$. By the genericity of the choice 
 $P$ in $S^d_x$, this fact would imply that 
 $$
  \overline T= H^0 ({\cal O}_{\ell}(d)(-x)),
 $$ 
 thus leading to a contradiction.
  
 Assume now that ${\ell_Q} \not ={\ell}_{Q'}$, 
 for generic $Q,Q'\in S^{d-1}$.
 Since $T$ contains a codimension two subspace of 
 $I_{\ell_Q} (d)$, from the exact sequence 
 $$
  0\rightarrow I_{\ell_Q} (d)\cap I_{{\ell}_{Q'}} (d)\rightarrow 
  I_{\ell_Q}(d)\oplus I_{\ell_{Q'}}(d)\rightarrow S^d_x\rightarrow 0,
 $$
 and the fact that $T\subset S^d_x$ has codimension $4$, 
 it follows that $T\supset I_{\ell_Q} (d)\cap I_{{\ell}_{Q'}} (d)$.
 Let $\mathbf P^2_{Q,Q'}$ be the span of ${\ell_Q}$ and
 $\ell_{Q'}$: we study the variation of this plane with $Q'$. 
 If for generic $Q'_1\not = Q'$ the intersection 
 $\mathbf P^2_{Q,Q'}\cap\mathbf P^2_{Q,Q'_1}$ is equal to the line 
 ${\ell}_Q$, then $T\supset I_{{\ell}_Q} (d)$ and we are done. 
 If otherwise $\mathbf P^2_{Q,Q'}=\mathbf P^2_{Q,Q'_1}= \mathbf P^2$,
 then it is immediate to see that $T$ contains
 $$
  \lbrace F\in S^d_x: F_{|\mathbf P^2} {\it\ is\ singular\ at\ the\ point} 
  \ x\rbrace,
 $$ 
 because ${\ell_Q}$ and $\ell_{Q'}$ will vary in this plane.
 But this is absurd since $T\subset S^d_x$ is of codimension $4$.
\end{proof}
From the previous lemma and the vanishing of the composite
map
$$
 H^0 (\wedge^{n-5-k} T{{\cal X}_U}_{|{X_F}}(n-5-k))\otimes 
 H^0 (\wedge^{4} T{{\cal X}_U}_{|{X_F}}(2))
 \rightarrow H^0 (K_{\tilde Y_F})
$$
we have that any codimension four subspace 
$T\subset T^{vert}_{{\cal X}_U,(y,F)}$ 
containing $T^{vert}_{{\cal Y},(y,F)}$ contains a hyperplane
$H_{\ell _{(y,F)}}$ of $I_{\ell _{(y,F)}}(d)$, where
$\ell_{(y,F)}$ is a line  through $y$. Note that {\it a priori} the
hyperplane $H_{\ell_{(y,F)}}$ could vary with $T$. 
We have then to verify 
that $T_{{\cal Y},(y,F)}$ is forced to contain one of those
$H_{\ell _{(y,F)}}$.

\bigskip
\noindent
{\it Proof of Proposition \ref{ideals}.} 
Remark that a codimension four subspace 
$T\subset T^{vert}_{{\cal X}_U,(y,F)}$ containing $T^{vert}_{{\cal Y},(y,F)}$ 
cannot contain two hyperplanes $H$ and $H'$ in the ideals 
$I_{\ell}(d)$ and $I_{\ell'}(d)$ of different lines ${\ell}\not ={\ell'}$. 
Indeed, if this was the case, by the surjectivity of
$$
 I_{\ell}(d)\oplus I_{\ell'}(d)\twoheadrightarrow S^d_y,
$$
then $T$ would contain at least a codimension 2 subspace of $S^d_y$,
thus violating the hypothesis on its codimension. 
Set $V:=S^d_y$, and $V_0:=T^{vert}_{{\cal Y},(y,F)}$. 
Denote by $\ell_T$ the unique line such that $T$ contains an
hyperplane in its ideal. Then,
by the above, we have a morphism from $G':=Grass (V/V_0,4)$, 
the Grassmannian of codimension 4 subspaces of $V/V_0$, to
$G:=Grass(1,{\mathbf P}^n)$, the Grassmannian of lines in ${\mathbf P}^n$:
$$
 \varphi\ :\ G'\to G;\ T\longmapsto \ell_T.
$$
Suppose that $\varphi$ is not constant map.
Let $\ell$ be in its image, and $T\in \varphi^{-1}(\ell)$. Then 
it is easy to construct from such a $T$ a positive dimensional 
family of codimension 4 subspaces containing $V_0$ and a hyperplane
in $I_{\ell}(d)$. 
Thus $\varphi$  has positive dimensional fibers and 
we are done, since in this case
the ample line bundle 
$\varphi^* {\cal O}_G (1)= {\cal O}_G' (s),\ s>0$, 
would have zero intersection with the curves in the fiber, which is 
absurd. 
Now, let $\ell$ be the unique line in the image of $\varphi$. 
A dimension count shows that if $V_0 \cap I_{\ell} (d)$ had 
codimension $\geq 2$, then it would exist a $T\in Grass(V/V_0, 4)$ 
with $codim_{I_{\ell} (d)} T\cap I_{\ell} (d)\geq 2$, thus contradicting
Lemma \ref{BL}.

\hfill $\Box$

%
\subsection{The vertical contact distribution}
%
We now want to use Proposition $3.1$ to construct a 
(well defined) distribution in
$T^{vert}_{\cal Y}$, and show its integrability.

From (\ref{ideals}) we know that $T_{{\cal Y},(y,F)}$
contains at least a hyperplane $H:=H_{\ell _{(y,F)}}$ in 
$I:=I_{\ell _{(y,F)}}(d)$. 
Remark that $T_{{\cal Y},(y,F)}$ cannot contain two different 
hyperplanes $H\not =H'$ in different ideals $I\not =I'$,
otherwise it would contain a codimension 
two subspace of $S^d_y$. But this is absurd, since by Lemma \ref{GL}
$$
 codim_{\ T^{vert}_{{\cal X},(y,F)}}\ T^{vert}_{{\cal Y},(y,F)}
 \ =\ codim_{\cal X}{\cal Y}\ \geq 4.
$$
Hence the line $\ell _{(y,F)}$ is unique and we have a well
defined map
\begin{eqnarray}
 \phi : &{\cal Y}&\longrightarrow\ \  G(1,n)\\
 &(y,F)& \longmapsto\ \  \ell _{(y,F)}\nonumber
\end{eqnarray}
If, at a generic point $(y,F)$, 
$T_{{\cal Y},(y,F)}$ contains the hyperplane 
$H_{\ell _{(y,F)}}$, but not the
whole ideal $I_{\ell _{(y,F)}}(d)$, we get a well defined distribution 
${\cal H}\subset T^{vert}_{\cal Y}$, whose fiber at a point $(y,F)$
is given by $H_{\ell _{(y,F)}}$. We will call ${\cal H}$ the 
{\it vertical contact distribution}.

If at a generic point 
$T_{{\cal Y},(y,F)}$ contains $I_{\ell _{(y,F)}}(d)$, then one can consider 
the distribution ${\cal I}\subset T^{vert}_{\cal Y}$ fiberwise defined
by $I_{\ell _{(y,F)}}(d)$. This case is easier and is actually
the case considered in [V3]. It will be briefly
treated at the end of this section.

In the former case, as in [V3] and by simply adapting the arguments given 
there to our situation, we now want to show the following natural fact: 
if we move infinitesimally in the directions parametrized by  
$H_{\ell _{(y,F)}}\subset I_{{\ell}_{(y,F)}}(d)$, then the 
line $\ell_{(y,F)}$ remains fixed. 
The integrability of $\cal H$ will then
immediately follow.  
\begin{lem}\label{integrability}
 
\begin{description} 
  \item{(i)} The differential $\phi_*$ at the point $(y,F)$ 
             vanishes on $H_{{\ell}_{(y,F)}}$.
  \item{(ii)} The vertical contact distribution 
  ${\cal H}\subset T^{vert}_{\cal Y}$ is integrable.
 \end{description} 
\end{lem}

\begin{proof}
 Since the distribution $T^{vert}_{{\cal Y},(y,F)}=ker\ p_*$ 
 is integrable, the brackets induce a map 
 $$
  {\Psi}:\bigwedge^2 {\cal H}\rightarrow
  T^{vert}_{\cal Y}/{\cal H}\subset 
  {T^{vert}_{\cal X}}_{|{\cal Y}}/{\cal H},
 $$
 which is given at the point $(y,F)$ by 
 $$
  \psi : \wedge^2 H_{\ell _{(y,F)}}
  \rightarrow T^{vert}_{{\cal Y},(y,F)}\ mod\ H_{\ell _{(y,F)}}
  \subset S^d_y\ mod\ H_{\ell _{(y,F)}}.
 $$
 Since we are supposing that $T^{vert}_{{\cal Y},(y,F)}$ contains $H_{(y,F)}$ 
 but not the whole ideal $I_{{\ell}_{(y,F)}}(d)$, there is a canonical 
 isomorphism 
 $$
  T^{vert}_{{\cal Y},(y,F)}\ mod\ H_{\ell _{(y,F)}}\cong 
  T^{vert}_{{\cal Y},(y,F)}\ mod \ I_{{\ell}_{(y,F)}}(d),
 $$
 and hence $\psi$ identifies with a map 
 $$
  \wedge^2 H_{\ell _{(y,F)}}
  \rightarrow H^0 ({\cal O}_{\ell} (d) (-y))
 $$
 which we also denote by $\psi$.
 To prove
 the integrability of $\cal H$, by Frobenius' theorem it will
 suffice to show that $\Psi$ is zero.
 In what follows we will denote 
 ${{\ell}_{(y,F)}}$ and $H_{(y,F)}$ respectively by $\ell$ and $H$.
 Now, choose coordinates on $\mathbf P^n$ such that 
 $\ell= \lbrace X_2=\ldots =X_n =0\rbrace$ and $y=[1,0,\ldots ,0]$.
 Recall that $H^0 (N_{\ell/{\mathbf P^n}}(-1))$ identifies naturally 
 with the set of $(-1)$-graded homomorphisms from 
 $\bigoplus_d I_{\ell}(d)/I_{\ell}^2(d)$ to 
 $\bigoplus_d S^d /I_{\ell}(d)$. Hence there is a natural bilinear
 map, denoted by $(a,b)\mapsto a\cdot b$ :
 $$
  I_{\ell} (d)\otimes H^0 (N_{\ell/{\mathbf P^n}}(-y))
  \rightarrow H^0 ({\cal O}_{{\ell}} (d) (-y)),
 $$ 
 which is explicitely given by 
 $$ 
  P\cdot (X_1 \sum_{i=2}^n b_i {\partial\over{\partial X_i}}) = 
  \sum_{i=2}^n b_i X_1({\partial P\over{\partial X_i}})_{|\ell} 
  \in H^0 ({\cal O}_{{\ell}} (d) (-y)).
 $$
 Remark that, since 
 $y\in \ell$ and $T^{vert}_{{\cal Y},(y,F)}\subset S^d_y$, 
 any deformation of $\ell$ belonging 
 to $\phi_* (T^{vert}_{{\cal Y},(y,F)})$ passes through $y$, i.e. 
 $\phi_* (T^{vert}_{{\cal Y},(y,F)})
 \subset H^0 (N_{\ell/{\mathbf P^n}}(-y))$.
 A verification in local coordinates shows that
 \begin{eqnarray}\label{star}
  \psi (A\wedge B)= A\cdot \phi_* (B)-B\cdot \phi_* (A),\ A,B\in H.
 \end{eqnarray}
 Note that 
 $$
  (QX_iX_j)\cdot (\sum_2^n b_l{\partial\over{\partial X_l}})= 
  \sum_2^n b_l({\partial QX_iX_j\over{\partial X_l}})_{|\ell}=0,
 $$
 for every $Q\in S^{d-2}$ and $i,j\geq 2$, 
 and therefore $\phi_*(A)\cdot B=0$, 
 for every $A\in H\cap I_{\ell}^2 (d),\ B\in H$.
 We first show that $\phi_*$ vanishes on $I_{\ell}^2 (d)$:
 if we had $\phi_*(A)\not =0$ with $A\in H\cap I_{\ell}^2 (d)$, then 
 $T^{vert}_{{\cal Y},(y,F)} mod\ I_{\ell} (d)$  would contain the elements
 $B\cdot\phi_*(A)$ for any $B\in H$, hence at least a hyperplane 
 of $H^0 ({\cal O}_{{\ell}} (d) (-y))$.
 Thus $\phi_*$ vanishes on $H\cap I^2_{\ell}(d)$, giving a map
 $H/H\cap I^2_{\ell}(d)\rightarrow H^0 (N_{\ell/{\mathbf P^n}}(-y))$, 
 which we still call $\phi_*$.\\
 Identify $H^0 ({\cal O}_{{\ell}} (d) (-y))$ 
 with $H^0 ({\cal O}_{{\ell}} (d-1))$, and recall again the natural
 isomorphism 
 $$ 
  I_{\ell}(d)/I^2_{\ell}(d)\cong 
  H^0 ({\cal O}_{{\ell}} (d-1))\otimes 
  H^0 (N_{\ell/{\mathbf P^n}}(-y))^*.
 $$
 Then $H/H\cap I^2_{\ell}(d)$ corresponds to a subspace
 $$
  {\overline H}\subset H^0 ({\cal O}_{{\ell}} (d-1))
  \otimes H^0 (N_{\ell/{\mathbf P^n}}(-y))^*,
 $$
 with $codim\ {\overline H}\leq 1$, and the dot map is 
 simply given by the contraction 
 $$
 <\cdot,\cdot>\ : {\overline H}\otimes H^0 (N_{\ell/{\mathbf P^n}}(-y))
  \rightarrow H^0 ({\cal O}_{{\ell}} (d)(-y))
 $$
 Hence, by (\ref{star}) the map
 $
 {\overline \psi} : \bigwedge^2 H/H\cap I^2_{\ell}(d)
 \rightarrow H^0 ({\cal O}_{{\ell}} (d)(-y))
 $
 identifies with
 \begin{eqnarray}
  && \bigwedge^2 {\overline H} \longrightarrow 
  H^0 ({\cal O}_{{\ell}} (d)(-y)) \\
  && A\wedge B \mapsto <\!\!A,\phi_*(B)\!\!>-<\!\!B,\phi_*(A)\!\!>. \nonumber
 \end{eqnarray}
 To conclude we need the following linear algebra result:
 \begin{lem}\label{lin.alg.}
  Let $W$ and $K$ be two vector spaces, 
  $\overline H$ a codimension $1$ subspace of in $W\otimes K^*$ 
  and $\phi_* :{\overline H}\rightarrow K$ a linear map.
  If $\phi_* \not = 0$, then the image of the map 
  \begin{eqnarray}
   {\overline \psi} &:& \bigwedge^2 {\overline H} 
   \longrightarrow W \nonumber \\
   && A\wedge B \mapsto <A,\phi_*(B)>-<B,\phi_*(A)> \nonumber
  \end{eqnarray}
  contains at least a codimension $2$ subspace of $W$. 
 \end{lem}
 \begin{proof}
  Let $J={\overline \psi}(\wedge^2 {\overline H})$. Pick a complementary 
  subspace $J^{\bot}$ to $J$ in $W$ and a basis 
  $\left\{ w_{j}\right\} $ for $W$ which is compatible with the decomposition
  $$
   W=J\oplus J^{\bot }.
  $$
  Let $\left\{ k_{i}\right\} $ be a basis of $K$ and $\{k^{*}_i\}$ the
  dual one.
  Pick a complementary space ${\overline H}^{\bot }$ to $\overline H$ 
  in $W\otimes K^*$ which will be generated by a monomial 
  $
   w_{j_0}\otimes k^{*}_{i_0},
  $ 
  and extend $\phi_*$ to the whole $W\otimes K^*$ by setting
  $$
   \phi_*(w_{j_0}\otimes k^{*}_{i_0})=0.
  $$
  The map ${\overline \psi}$ extends naturally to $\wedge^2 (W\otimes K^*)$.
  Since 
  $$ 
  \dim \left( \psi \left( {\overline H}^{\bot }
  \otimes \left( W\otimes K^*\right) \right)
  \right) \leq 1,
  $$
  we are reduced to proving that if the extended map $\phi_*:W\to K$ 
  is not zero, then the codimension of the image of
  $$
   \overline\psi\ :\ \bigwedge^2 W\otimes K^*\to W
  $$
  is at least $1$.
  This has already been checked in [V3], Lemma $3$, and so we are done.
 \end{proof}
 Take $W:=H^0 ({\cal O}_{{\ell}} (d)(-y))$ and 
 $K:=H^0 (N_{\ell/{\mathbf P^n}}(-y))$, and apply Lemma $\ref{lin.alg.}$
 to our situation. If we had $\phi_*\not =0$, then the image of the 
 map ${\overline \psi}$ would contain at least a codimension $2$ subspace of  
 $H^0 ({\cal O}_{{\ell}} (d)(-y))$. But the image of ${\overline \psi}$
 is contained in $T^{vert}_{{\cal Y},(y,F)}\ mod\ I_{\ell}(d)$, 
 and $T^{vert}_{{\cal Y},(y,F)}$ contains a hyperplane in $I_{\ell}(d)$.
 Hence the codimension of 
 $T^{vert}_{\cal Y}$ in $T^{vert}_{\cal X}$ would be at most $3$, which is 
 in contradiction with Lemma \ref{pos}.
 Thus $\phi_*=0$, hence $\psi$ is zero and so is $\Psi$. 
 By Frobenius' theorem the distribution $\cal H$ is integrable.
\end{proof}

%
%
\subsection{Proof of Proposition A}
%
%

We can now prove Proposition A, i.e. we show that the line
${\ell}_{(y,F)}$  defined by the ideal $I_{{\ell}_{(y,F)}}(d)$ 
is such that
$$ 
 X_F\cap {\ell}_{(y,F)}=d\cdot y.
$$

From (\ref{integrability}) we know that $\cal H$ is integrable and 
$\phi$ is constant along the leaves of the corresponding foliation. 
Therefore the line ${\ell}_{(y,F)}$ is fixed along the leaf, and 
because its tangent space is contained at each point in 
$I_{{\ell}_{(y,F)}}(d)$,
it follows that the restriction $G_{|{\ell}_{(y,F)}}$ is 
constant, for any polynomial $G$ belonging to the leaf through $(y,F)$.
This means that the leaf is locally of the form 
$y\times F+W_{(y,F)}$, where $W_{(y,F)}\subset I_{{\ell}_{(y,F)}}(d)$ is 
a germ of complex hypersurface. Then consider the
following diagram
\begin{equation}\label{codim}
\xymatrix{ 
 &0\ar[d]
&0\ar[d]
\\
0\ar[r]
&H_{{\ell}_{(y,F)}}\ar[d]\ar[r]
&I_{{\ell}_{(y,F)}}(d)\ar[d]
\\
0\ar[r]
&T^{vert}_{{\cal Y},(y,F)}\ar[d]\ar[r]
&T^{vert}_{{\cal X},(y,F)}=S^d_y\ar[d]
\\
 &H^0 ({\cal O}_{{\ell}} (d)(-y))\ar@{=}[r]
&H^0 ({\cal O}_{{\ell}} (d)(-y))\ar[d]
\\
 &
 &0
}\end{equation}
By Lemma \ref{GL},
$$
 codim_{\ T^{vert}_{{\cal X},(y,F)}}\ T^{vert}_{{\cal Y},(y,F)}
 =codim _{\ \cal X}\ {\cal Y}=n-k-1.
$$
Then, by (\ref{codim}) the image 
$$
 Im:=\ Im\ (T^{vert}_{{\cal Y},(y,F)}
 \rightarrow H^0 ({\cal O}_{{\ell}} (d)(-y)))
$$
has codimension $(n-k-1)-1=n-k-2$, and therefore, since $d=2n-2-k$ 
\begin{eqnarray}\label{dim.}
 dim\ Im\ =(2n-2-k)-(n-k-2)=n.
\end{eqnarray}
At the same time, again by Lemma \ref{GL}, 
$T^{vert}_{{\cal Y},(y,F)}$ contains
$S^1_y \cdot J_F^{d-1}$ and $F$ itself. 
Take coordinates $X_0,\ldots ,X_n$ on
$\mathbf P^n$ such that $y=[1,0,\ldots ,0]$, and
${\ell}:={\ell}_{(y,F)}=\lbrace X_2=\ldots =X_n=0\rbrace$.
Since $\phi$ is constant along the leaves of the foliation, 
we can generically choose a polynomial $G$ in
$F+ W_{(y,F)}$, so that the $(n-1)$-elements
$X_1{\partial G\over{\partial X_i}},\ i\geq 2$, 
are generic in a hypersurface.
Consider the subspace 
$$
 K:= <G_{|{\ell}},
 X_1({\partial G\over{\partial X_0}})_{|{\ell}},
 X_1({\partial G\over{\partial X_1}})_{|{\ell}}>\ 
 \subset H^0 ({\ell}, {\cal O}_{{\ell}}(d)\otimes
 {\cal I}_y),
$$
which is uniquely determined by $F_{|{\ell}}$ and hence is
constant along the leaf.  Its codimension
in $H^0 ({\ell}, {\cal O}_{\ell}(d)\otimes {\cal I}_y)$ 
is at least $d-3\geq n$ 
(since, by hypothesis, $k\leq n-5$ so that $d=2n-2-k \geq n+3$). 
Since we know that along the leaf, $G$ moves freely in the complex 
hypersurface $W_{(y,F)}$, the polynomials 
$X_1({\partial G\over{\partial X_i}})_{|{\ell}}$ are generic in a 
codimension $1$ subspace of 
$H^0 ({\ell}, {\cal O}_{\ell}(d)\otimes {\cal I}_y)$, and it then
follows that they will be generically independant modulo $K$.
From (\ref{dim.}) we thus get that 
\begin{eqnarray}\label{dim.2}
 dim\ K\leq 1
\end{eqnarray}
and so Proposition A is proved, since by (\ref{dim.2}) $F_{|{\ell}}$ has
to be of the form $\alpha X_1^d$.

If, for generic $(y,F)$, $T^{vert}_{{\cal Y},(y,F)}$ contains 
the whole ideal $I_{{\ell}_{(y,F)}}(d)$, then consider the distribution
${\cal I}\subset T^{vert}_{\cal X}$  
pointwise defined by $I_{\ell _{(y,F)}}(d)$.
Arguing as we did before, we get 
$$
  dim\ Im\ (T^{vert}_{{\cal Y},(y,F)}
  \rightarrow H^0 ({\cal O}_{{\ell}} (d)(-y)))=n-1,
$$
thus deducing 
$$
 dim\ K =0.
$$
Then the polynomial $F$ belongs to $I_{{\ell}_{(y,F)}}(d)$, 
and the theorem is true in this case, 
i.e. $Y_F$ is a component of the subvariety
of $X_F$ covered by lines.\hfill $\Box$
%
\section {The geometry of $\Delta_{d,F}$}
%
Let $X_F\subset \mathbf P^n$ be a general 
hypersurface of degree $d={2n-2-k},\ 1\leq k\leq n-5$, and $Y_F\subset X_F$ a
$k$-dimensional subvariety  whose desingularization 
$\tilde Y$ has $h^0 (\tilde Y,K_{\tilde Y})=0$. 
Then, by Proposition A, we know that
$Y_F$ has to be contained in $\Delta_{d,F}\subset X_F$, 
the $(k+1)$-dimensional subvariety
of points in $X_F$ through which there is a $d$-osculating line.
To prove Proposition B and hence our theorem, we have then to show
that the only subvariety of dimension $k$ of $\Delta_{d,F}$ with geometric
genus zero is the subvariety covered by the lines in $X_F$.

%
\subsection{A desingularization of $\Delta_{d,F}$}
%

We start by giving an explicit description of a desingularization
$\tilde\Delta_{d,F}$ of $\Delta_{d,F}$ in terms of the zero locus 
of a section of a vector bundle.
This fact will allow us to calculate, by adjunction,
the canonical bundle of $\tilde\Delta_{d,F}$ and see that it is very ample.

Let $G:= Gr(1,n)$ be the Grassmannian of lines in $\mathbf P^n$. 
Let ${\cal O}_G (1)$ be the line bundle on $G$ which gives the Pl\"ucker 
embedding, so that we have
$H^0 ({\cal O}_G (1))= \bigwedge^2 S^1$.
Let ${\cal P}\subset {\mathbf P^n} \times G$ be the incidence variety 
$\lbrace (x,[\ell]) :\ x\in \ell\rbrace$ with projections
\begin{equation}
\xymatrix{
{\cal P}\ar[r]^{p}\ar[d]^{q}
& \mathbf P^n\\
G
& \\
}\end{equation}
and $H:=p^* {\cal O}_{\mathbf P^n}(1)$, $L:=q^* {\cal O}_G (1)$ the 
line bundles generating the Picard group of $\cal P$.

Define 
$$
 \tilde \Delta_r:=\{(x,[\ell],F)\ :\ \ell\cdot X_F\geq r\cdot x \}
 \subset {\cal P}\times S^r
 \subset \mathbf P^n\times G\times S^r,
$$
(by $\ell\cdot X_F\geq r\cdot x$ we mean that the line has a 
contact of order at least $r$ with $X_F$ at $x$)
and consider the various projections as illustrated in the following
commutative diagram:
\begin{equation}\label{proj}
\xymatrix{ 
\tilde \Delta_r\ar[r]^{\rho_r}\ar[d]^{\pi}
&{\cal X}\ar[d]^{t}\ar[r]^{s}
&S^r
\\
{\cal P}\ar[r]^{p}\ar[d]^{q}
&{\mathbf P}^n
&
\\
G &
&
}\end{equation}

Since the tangency of order at least $r$ imposes $r$ conditions, 
the fibres of the projection $\pi :\tilde\Delta_r\to{\cal P}$ are punctured vector 
spaces of dimension $N-r$. Hence $\tilde\Delta_r$ is smooth and irreducible
of dimension
$$
 N-r+2(n-1)+1.
$$
\begin{lem}\label{des}
(i) The projection $\rho_r: \tilde\Delta_r\to {\cal X}$
is surjective for $r\leq n$, and generically injective for $r>n$. 

(ii) The composite projection 
$
s\circ\rho_r\ :\tilde\Delta_r\to S^r
$
is surjective if $r\leq 2(n-1)$. In particular, in that case, 
its fiber $\tilde\Delta_{r,F}:=s\circ\rho_r^{-1}(F)$
is smooth for generic $F\in S^r$, and the composite projection
$t\circ \rho_r:\tilde\Delta_{r,F}\to {\mathbf P}^n$, 
gives a desingularization of  
$$
\Delta_{r,F}:=\{x\in X_F\ :\ \exists \ell\ s.t.\ \ell\cdot X_F\geq r\cdot x \}
$$ 
\end{lem}
\begin{proof}
 (i) Assume
$
x=\left[ 1,0,\ldots ,0\right] \in Proj\left( \Bbb{C}\left[ X_{0},\ldots
,X_{n}\right] \right) .
$
Then the assertion follows from the fact that the contact condition 
$\ell\cdot X_F\geq r$ for a line $\ell$ through $x$ with respect to
$$
 F=\sum\nolimits_{j=1}^{n}X_{0}^{d-j}F_{j}\left( X_{1},\ldots ,X_{n}\right)
$$
becomes
$$
 \left\{ F_{1}=\ldots =F_{r-1}=0\right\} \subseteq Proj\left( \Bbb{C}\left[
 X_{1},\ldots ,X_{n}\right] \right).
$$
 (ii) A dimension count shows that all hypersurfaces $X_{F}$ in 
${\mathbf P}^{n}$ of degree $d\leq 2n-2$
admit a point through which passes a line having contact with $X_F$ of 
maximal order.
\end{proof}

In what follows, by abuse of notation, we will identify $\tilde\Delta_{r,F}$
to its image, $\pi (\tilde\Delta_{r,F})$, in ${\cal P}$.
We will show that $\tilde\Delta_{r,F}$ can be seen as the zero locus
of a global section of a vector bundle over ${\cal P}$. This will enable us
to compute its canonical bundle.

Let ${\cal E}_d$ be the $d^{th}$-symmetric power of the dual of the 
tautological subbundle on $G$, 
and recall that, by definition, its fibre at a point $[\ell]$ 
is then given by $H^0 ({\ell},{\cal O}_{\ell} (d))$,
and its first Chern class is
$$
 c_1({\cal E}_d)={\cal O}_{G} ({d(d+1)\over 2}).
$$
 Let
${\cal L}_d :=dL-dH$ be the rank $1$ subbundle of $q^* {\cal E}_d$. Note 
that its fibre ${\cal L}_{d,(x,[\ell])}$ 
is equal to the space of degree $d$ 
homogeneous polynomials on $\ell$ vanishing to the order $d$ at $x$. 
Finally, let ${\cal F}_d$ be the quotient 
\begin{eqnarray}\label{L,E,F}
 0\rightarrow {\cal L}_d\rightarrow q^* {\cal E}_d 
 \rightarrow {\cal F}_d\rightarrow 0.
\end{eqnarray}
It is possible to associate to every $F\in S^d$ a section 
$\sigma_F\in H^0 (G,{\cal E}_d)$, whose value at a point 
$[\ell]$ is exactly the polynomial $F_{|\ell}$. 
We will denote by ${\overline {\sigma}}_F$ 
the induced section in $H^0 ({\cal P},{\cal F}_d)$.
Then $V({\overline {\sigma}}_F)$ is equal to 
$\tilde\Delta_{r,F}$, and, as checked in Lemma \ref{des}, 
for generic $F$,
$\tilde\Delta_{r,F}=V({\overline {\sigma}}_F)
{\buildrel p\over \longrightarrow} \Delta_{r,F}$ 
is a desingularization. It is computed in [V3] that 
\begin{equation}\label{V}
 K_{\cal P}=-2H-nL,
\end{equation} 
and then
\begin{equation}\label{det}
  det\ {\cal F}_d= det\ q^* {\cal E}_d -dL+dH={d(d-1)\over{2}} L+dH.
\end{equation} 
Hence, by adjunction, 
\begin{equation}\label{K}
 K_{{\tilde\Delta_{d,F}}}= 
 (K_{\cal P}\otimes det\ {\cal F}_d)_{|{\tilde \Delta_{d,F}}}= 
 [(d-2)H+ ({d(d-1)\over 2}-n)L]_{|{\tilde\Delta_{d,F}}},
\end{equation} 
so $K_{\tilde\Delta_{d,F}}$ is very ample under our numerical
hypothesis $d=2n-2-k,\ 1\leq k\leq n-5$.

%
\subsection{Proof of Proposition B}
%

Let $\Delta_d\subset \mathbf P^n\times S^d$ be the family of the 
$\Delta_{d,F}$'s,
and $\tilde\Delta_d\subset {\cal P}\times S^d$ the family of the 
desingularizations. Let ${\cal Y}\subset \tilde\Delta_d$ be a subscheme
of relative dimension $k$, invariant under the action of
$GL(n+1)$, and $\tilde {\cal Y}\rightarrow {\cal Y}$ a desingularization.
Assume $h^0 ({\tilde Y}_F,K_{{\tilde Y}_F})=0$.
Recall the isomorphisms 
\begin{eqnarray}
 {{T{\tilde \Delta_d}}}_{|{{\tilde\Delta}_{d,F}}}
 \otimes K_{\tilde\Delta_{d,F}}&\cong&
 {{\Omega}^{N+k}_{{\tilde \Delta_d}}}_{|\tilde \Delta_{d,F}}\\ 
 {\Omega ^{N+k}_{{\tilde{\cal Y}}}}_{|{\tilde Y_F}}
 &\cong& K_{\tilde Y_F}
\end{eqnarray}
and consider the natural map 
\begin{eqnarray}
 {{T{\tilde \Delta_d}}}_{|{{\tilde\Delta}_{d,F}}}
 \otimes K_{\tilde\Delta_{d,F}}\cong
 {{\Omega}^{N+k}_{{\tilde \Delta_d}}}_{|\tilde \Delta_{d,F}}\to
 {\Omega ^{N+k}_{{\tilde{\cal Y}}}}_{|{\tilde Y_F}}
 &\cong& K_{\tilde Y_F}.
\end{eqnarray}
Then, by assumption, the induced map in cohomology 
\begin{equation}\label{S-restriction}
 H^0 ({{T{\tilde \Delta_d}}}_{|{{\tilde\Delta}_{d,F}}}
 \otimes K_{\tilde\Delta_{d,F}})
 \rightarrow H^0 (K_{{\tilde Y}_F})
\end{equation}
is zero. Let $T^{vert}_{\tilde \Delta_d}$ be the sheaf defined by
$$
 0\rightarrow T^{vert}_{\tilde \Delta_d}
 \rightarrow T{\tilde \Delta_d}
 \buildrel{\pi _*}\over{\longrightarrow} T{\cal P}
 \rightarrow 0.
$$
Using the positivity result proved in Lemma \ref{pos}, (ii), we will
construct a sub-bundle of
$
(T^{vert}_{\tilde \Delta_d})_{|{\tilde \Delta_{d,F}}}
\otimes K_{\tilde\Delta_{d,F}},
$ 
generated by its global sections. This will allow us to 
show that any point $(y,[\ell],F)\in {\cal Y}$ is 
such that $y\in\ell\subset X_F$. Comparing the dimension, we will
thus obtain that $Y_F$ has to be a component of the subvariety 
of lines in $X_F$.

From (\ref{S-restriction}) we see that, at a smooth point 
$(y,[\ell],F)\in {\cal Y}\subset {\cal P}\times S^d$, the tangent space 
$T_{{\cal Y},(y,[\ell],F)}$ is in the base locus of 
$H^0 ({{T{\tilde \Delta_d}}}_{|{{\tilde\Delta}_{d,F}}}
 \otimes K_{\tilde\Delta_{d,F}})$, 
considered
as the space of sections of a line bundle on the Grassmannian of
hyperplanes in ${{T{\tilde \Delta_d}}}_{|{{\tilde\Delta}_{d,F}}}$.
Consider the vector bundle $M^d_G$ on $G:= Gr(1,n)$ 
defined by the short exact sequence:
$$
 0\rightarrow M^d_G\rightarrow S^d \otimes {\cal O}_G 
 \rightarrow {\cal E}_d \rightarrow 0.
$$
Notice that the fiber of $M^d_G$ at a point 
$[\ell]$ is equal to $I_{\ell}(d)$, and recall that, 
by Proposition \ref{pos}, (ii), $M^d_G\otimes {\cal O}_G(1)$ is generated by
its global sections.

Then  it follows that the vector bundle 
$q^* M^d_G \otimes det {\cal F}_d\otimes K_{\cal P}$, that, by 
(\ref{V}), (\ref{det}) and (\ref{K}) can be written as
\begin{eqnarray}
 q^* M^d_G \otimes det {\cal F}_d\otimes K_{\cal P}
 \!\!\! &=&\!\!\!\! q^* (M^d_G) \otimes {\cal O}_{\cal P}((d-2)H+ 
 ({d(d-1)\over 2}-n)L)\nonumber\\
 &=& \!\!\!\! q^* (M^d_G (1)) \otimes {\cal O}_{\cal P}((d-2)H+ 
 ({d(d-1)\over 2}-n-1)L),\nonumber
\end{eqnarray}
is generated by its global sections, and so is its restriction 
to ${\tilde\Delta}_{d,F}$, 
i.e. the sheaf 
$$
 {q^*{M^d_G}}_{|{{\tilde\Delta}_{d,F}}}
 \otimes K_{\tilde\Delta_{d,F}}
$$
is generated by its global sections. 

Let ${\cal N}_d$ be the vector bundle over $\cal P$ 
defined by the exact sequence
\begin{eqnarray}\label{N_d}
 0\rightarrow {\cal N}_d\rightarrow S^d \otimes {\cal O}_{\cal P}
 \rightarrow {\cal F}_d\rightarrow 0.
\end{eqnarray}
We have  
$$
 0\rightarrow {{\cal N}_d}_{|{\tilde \Delta_{d,F}}}
 \rightarrow 
 {{T{\tilde \Delta_d}}}_{|{{\tilde\Delta}_{d,F}}}
 \buildrel{\pi _*}\over{\longrightarrow} T{\cal P}_{|{\tilde \Delta_{d,F}}}
 \rightarrow 0,
$$
where 
${\cal S}\subset {\cal P}\times S^d \buildrel{\pi}\over{\rightarrow} {\cal P}$ 
is the projection  on the first component, 
i.e. ${{\cal N}_d}_{|{{\tilde\Delta}_{d,F}}}$ is
the vertical component of 
${{T{\tilde \Delta_d}}}_{|{{\tilde\Delta}_{d,F}}}$
w.r.t. $\pi$.
Now consider the vector bundle ${\cal M}^d_G$ defined by the exact sequence
$$
 0\rightarrow {\cal M}^d_G\rightarrow S^d\otimes {\cal O}_{\cal P}
 {\buildrel {ev}\over \longrightarrow}
 q^* {\cal E}_d\rightarrow 0,
$$
whose fiber at a point $(y,F,[\ell])$ 
is equal to $I_{\ell}(d)$.
From (\ref{L,E,F}) and the definition of ${\cal N}_d$ we also obtain that
$$
 0\rightarrow {\cal M}^d_G\rightarrow 
 {\cal N}_d\rightarrow {\cal L}_d\rightarrow 0.
$$
So, ${{\cal M}^d_G}\otimes K_{\tilde\Delta_{d,F}}$ is a subbundle
of ${{T{\tilde \Delta_d}}}^{vert}_{|{{\tilde\Delta}_{d,F}}}$.
Finally note that ${\cal M}^d_G= q^* M^d_G$,
hence 
$$
 {\cal M}^d_G\otimes K_{\tilde\Delta_{d,F}}
$$ 
is generated by its global sections. Using this property of the bundle
${{\cal M}^d_G}\otimes K_{\tilde\Delta_{d,F}}
\subset 
{{T{\tilde \Delta_d}}}^{vert}_{|{{\tilde\Delta}_{d,F}}}
 \otimes K_{\tilde\Delta_{d,F}}$
we are now able to conclude our proof.
\bigskip
\\
{\it Proof of Proposition B.}
 Let $H\subset {T_{\tilde \Delta_d}}_{,(x,\Delta,F)}$ 
 be a hyperplane contained 
 in the base locus of 
 $H^0 ({{T{\tilde \Delta_d}}}_{|{{\tilde\Delta}_{d,F}}}
 \otimes K_{\tilde\Delta_{d,F}})$, considered 
 as the space of sections of a line bundle
 on the Grassmannian of codimension $1$ subspaces of 
 ${{T{\tilde \Delta_d}}}_{|{{\tilde\Delta}_{d,F}}}$. 
 Then we must have
 \begin{eqnarray}\label{vert.bl}
  H^{vert}:=H\cap {\cal N}_{d,(x,[\ell])}={\cal M}^d_{G,(x,\ell)}.
 \end{eqnarray}
 Indeed, if  
 $\bar H:= H\cap {\cal M}^d_{G,(x,[\ell])}$ was strictly contained in 
 ${\cal M}^d_{G,(x,[\ell])}$,
 then consider the following, well defined, commutative diagram:
 \begin{equation}
  \xymatrix{
   H^0 ({\cal M}^d_G\otimes K_{\tilde\Delta_{d,F}})
   \ar[r]^{ev\ \ }\ar@^{(->}[d]\ar@_{>}[d]
   &({\cal M}^d_G\otimes K_{\tilde\Delta_{d,F}})_{(x,[\ell])}
   \ar[r]^{\ \ \ \ \ \ \ \ \ \ \ <\cdot,{\bar H}>}\ar@^{(->}[d]\ar@_{>}[d]
   &{\mathbb C}\ar@^{=}[d]\\
   H^0 ({{T{\tilde \Delta_d}}}_{|{{\tilde\Delta}_{d,F}}}
   \otimes K_{\tilde\Delta_{d,F}})\ar[r]^{ev\ \ }
   &({{T{\tilde \Delta_d}}}_{|{{\tilde\Delta}_{d,F}}}
   \otimes K_{\tilde\Delta_{d,F}})_{(x,[\ell])}
   \ar[r]^{\ \ \ \ \ \ \ \ \ \ \ \ \ \ \ <\cdot,H>}
   &{\mathbb C}\\
 }\end{equation}
 ($ev$ is the evaluation of the sections at the point $(x,[\ell])$, and
 $<\cdot,H>$ is the contraction defined by the hyperplane $H$).
 Since $H$ belongs to the base locus of 
 $H^0 ({{T{\tilde \Delta_d}}}_{|{{\tilde\Delta}_{d,F}}}
   \otimes K_{\tilde\Delta_{d,F}})$, then the composite map 
 $<\cdot,{H}>\circ\ ev$ is zero, and 
 so would be $<\cdot,{\bar H}>\circ\ ev$. But this is absurd,
 because 
 ${\cal M}^d_G\otimes K_{\tilde\Delta_{d,F}}$
 is generated by its global sections.

 Let then ${\cal Y}\subset {\Delta}_d$ be a subvariety, 
 which is stable under the action of $GL(n+1)$ and of relative 
 codimension $1$. Assume moreover 
 that the restriction map (\ref{S-restriction})
 $$
 H^0 ({{T{\tilde \Delta_d}}}_{|{{\tilde\Delta}_{d,F}}}
 \otimes K_{\tilde\Delta_{d,F}})
 \rightarrow H^0 (K_{{\tilde Y}_F})
 $$
 is zero.
 By (\ref{vert.bl}), $T^{vert}_{{\cal Y},(y,[\ell],F)}$ is equal to
 \begin{equation}\label{etoile}
 {\cal M}_{d, (y,[\ell])}=\lbrace G\in S^d : G_{|\ell}=0\rbrace.
 \end{equation} 
 On the other hand, by Lemma \ref{GL}, (ii), 
 $T^{vert}_{{\cal Y},(y,[\ell],F)}$ contains  
 $F$ itself.
 So by (\ref{etoile}) we have that $F_{|\ell}=0$ 
 for every point $(y,[\ell])\in Y_F$, 
 i.e. $Y_F$ is a component of the subvariety covered
 by the lines contained in $X_F$.  
 \hfill $\Box$

\begin{rmk}
If $k>1$, the $k$-dimensional subvariety covered by the lines of 
the general hypersurface of degree $d=2n-2-k$ is
irreducible (see [DM]), so in this case $Y_F$ has to coincide with it.
\end{rmk}

\noindent Gianluca PACIENZA
\\Institut de Math\'ematiques de Jussieu
\\Universit\'e Pierre et Marie Curie
\\2, Place Jussieu, F-75252 Paris CEDEX 05 - FRANCE
\\e-mail: pacienza@math.jussieu.fr 
\end{document}